\documentclass[amstex,12pt,russian,amssymb]{article}

\usepackage{mathtext}
\usepackage[cp1251]{inputenc}
\usepackage[T2A]{fontenc}
\usepackage[russian]{babel}
\usepackage[dvips]{graphicx}
\usepackage{amsmath}
\usepackage{amssymb}
\usepackage{amsxtra}
\usepackage{latexsym}
\usepackage{ifthen}

\begin{document}

\newcounter{lemma}
\newcommand{\lemma}{\par \refstepcounter{lemma}%
{\bf Лема \arabic{lemma}.}}

\newcounter{corollary}
\newcommand{\corollary}{\par \refstepcounter{corollary}%
{\bf Наслідок \arabic{corollary}.}}

\newcounter{remark}
\newcommand{\remark}{\par \refstepcounter{remark}%
{\bf Зауваження \arabic{remark}.}}

\newcounter{theorem}
\newcommand{\theorem}{\par \refstepcounter{theorem}%
{\bf Теорема \arabic{theorem}.}}

\newcounter{proposition}
\newcommand{\proposition}{\par \refstepcounter{proposition}%
{\bf Твердження \arabic{proposition}.}}

\newcounter{example}
\newcommand{\example}{\par \refstepcounter{example}%
{\bf Приклад \arabic{example}.}}

\renewcommand{\refname}{\centerline{\bf Список літератури}}

\renewcommand{\figurename}{Мал.}

\newcommand{\proof}{{\it Доведення.\,\,}}

\noindent УДК 517.5

{\bf А.П.~Довгопятый} (Житомирский государственный университет имени
Ивана Франко)

{\bf Е.А.~Севостьянов} (Житомирский государственный университет
имени Ивана Франко; Институт прикладной математики и механики НАН
Украины, г.~Славянск)

\medskip\medskip
{\bf О.П.~Довгопятий} (Житомирський державний університет імені
Івана Фран\-ка)

{\bf Є.О.~Севостьянов} (Житомирський державний університет імені
Івана Фран\-ка; Інститут прикладної математики і механіки НАН
України, м.~Слов'янськ)

\medskip\medskip
{\bf O.P.~Dovhopiatyi} (Zhytomyr Ivan Franko State University)

{\bf E.A.~Sevost'yanov} (Zhytomyr Ivan Franko State University;
Institute of Applied Ma\-the\-ma\-tics and Mechanics of NAS of
Ukraine, Slov'yans'k)

\medskip
{\bf О компактности классов решений задачи Дирихле с ограничениями
те\-о\-ре\-ти\-ко-множественного типа}

{\bf Про компактність класів розв'язків задачі Діріхле з обмеженнями
тео\-ре\-ти\-ко-множинного типу}

{\bf On compactness of classes of solutions of the Dirichlet problem
with restrictions of the theoretics-set type}

\medskip\medskip
Доказаны теоремы о компактных классах гомеоморфизмов с
гидродинамической нормировкой, являющихся решениями уравнения
Бельтрами, характеристики которых имеют компактный носитель и
удовлетворяют определённым ограничениям теоретико-множественного
типа. В качестве следствий, получены результаты о компактных классах
решений соответствующих задач Дирихле, рассматриваемых в некоторой
жордановой области.

\medskip\medskip
Доведено теореми про компактні класи гомеоморфізмів з
гідродинамічним нормуванням, які є розв'язками рівняння Бельтрамі,
характеристики яких мають компактний носій і задовольняють певні
обмеження теоретико-множинного типу. Як наслідок, отримано
результати про компактні класи розв'язків відповідних задач Діріхле,
які розглядаються в деякій жордановій області.

\medskip\medskip
We have proved theorems on compact classes of homeomorphisms with
hydrodynamic normalization that are solutions of the Beltrami
equation, whose characteristics are compactly supported and satisfy
certain constraints of the theoretical-set type. As a consequence,
we obtained results on compact classes of solutions of corresponding
Dirichlet problems considered in some Jordan domain.

\newpage
{\bf 1. Вступ.} Відносно нещодавно отримано результати щодо
компактності сімей розв'яз\-ків рівнянь Бельтрамі, а також
відповідної задачі Діріхле для нього (див., напр., \cite{Dyb} і
\cite{L$_2$}). Зокрема, отримано компактність цих розв'язків з
умовами нормування $f(0)=0,$ $f(1)=1$ і $f(\infty)=\infty,$
характеристики яких задовольняють обмеження інтегрального, або
теоретико-множинного типу. Окремо вивчалася задача Діріхле для
рівняння Бельтрамі. Зокрема, встановлено існування відкритих
дискретних розв'язків цієї задачі, а також теореми компактності їх
сімей в одиничному крузі (див., напр., \cite{Dyb} і \cite{Dyb$_2$}).
Метою даної публікації є отримання нових умов компактності класів
розв'язків рівняння Бельтрамі і задачі Діріхле в випадку, коли їх
характеристики задовольняють теоретико-множинні обмеження на
дилатації. Зокрема, доведено компактність класів розв'язків рівняння
Бельтрамі з так званим гідродинамічним нормуванням, тобто, коли ці
розв'язки ведуть себе близько до тотожних відображень в околі
нескінченно віддаленої точки. Аналогічні результати отримано для
розв'язків задачі Діріхле у довільній обмеженій жордановій області.

\medskip
Скрізь далі відображення $f:D\rightarrow {\Bbb C}$ області
$D\subset{\Bbb C}$ вважається таким, що {\it зберігає орієнтацію,}
зокрема, якщо $f$ -- гомеоморфізм і $z\in D$ -- яка-небудь його
точка диференційовності, то {\it якобіан} цього відображення в точці
$z$ додатній. Для комплекснозначної функції $f:D\rightarrow {\Bbb
C},$ заданій в області $D\subset {\Bbb C},$ що має частинні похідні
по $x$ і $y$ при майже всіх $z=x + iy,$ покладемо $f_{\overline{z}}
= \left(f_x + if_y\right)/2$ і $f_z = \left(f_x - if_y\right)/2.$
{\it Комплексною дилатацією} відображення $f$ в точці $z$
називається функція $\mu:D\rightarrow {\Bbb C},$ визначена рівністю
$\mu(z)=\mu_f(z)=f_{\overline{z}}/f_z$ при $f_z \ne 0$ і $\mu(z)=0$
в іншому випадку. {\it Мак\-си\-маль\-ною дилатацією} відображення
$f$ в точці $z$ називається наступна функція:
\begin{equation}\label{eq1}
K_{\mu}(z)=K_{\mu_f}(z)=\quad\frac{1+|\mu (z)|}{1-|\mu\,(z)|}\,.
\end{equation}
Якщо задана вимірна за Лебегом функція $\mu:D\rightarrow {\Bbb D},$
${\Bbb D}=\{z\in {\Bbb C}: |z|<1\},$ то не прив'язуючись до
якого-небудь відображення $f$ будемо називати величину, що
обчислюється за допомогою рівності~(\ref{eq1}), максимальною
дилатацією відповідної функції $\mu.$ Зауважимо, що якобіан
відображення $f$ в точці $z\in D$ можна обчислити за допомогою
рівності
$$J(z,
f)=|f_z|^2-|f_{\overline{z}}|^2\,,$$
що можна перевірити прямими обчисленнями. Неважко бачити, що
$K_{\mu_f}(z)=\frac{|f_z|+|f_{\overline{z}}|}{|f_z|-|f_{\overline{z}}|}$
у всіх точках $z\in D$ відображення $f,$ що має частинні похідні в
точці $z,$ де якобіан $J(z, f)$ не дорівнює нулю. {\it Рівнянням
Бельтрамі} будемо називати диференціальне рівняння виду
\begin{equation}\label{eq2}
f_{\overline{z}}=\mu(z)\cdot f_z\,,
\end{equation}
в якому $\mu=\mu(z)$ -- задана невідома функція. {\it Регулярним
розв'язком рівняння~(\ref{eq2})} в області $D\subset{\Bbb C}$ ми
будемо називати гомеоморфізм $f:D\rightarrow{\Bbb C}$ класу $W_{\rm
loc}^{1, 1}(D)$ такий, що $J(z, f)\ne 0$ при майже всіх $z\in D.$ У
подальшому, в розширеному просторі $\overline{{{\Bbb R}}^n}={{\Bbb
R}}^n\cup\{\infty\},$ $n\geqslant 2,$ використовується {\it сферична
(хордальна) метрика} $h(x,y)=|\pi(x)-\pi(y)|,$ де $\pi$ --
стереографічна проекція $\overline{{{\Bbb R}}^n}$  на сферу
$S^n(\frac{1}{2}e_{n+1},\frac{1}{2})$ в ${{\Bbb R}}^{n+1},$ а саме,
$$h(x,\infty)=\frac{1}{\sqrt{1+{|x|}^2}}\,,$$
\begin{equation}\label{eq3C}
\ \ h(x,y)=\frac{|x-y|}{\sqrt{1+{|x|}^2} \sqrt{1+{|y|}^2}}\,, \ \
x\ne \infty\ne y
\end{equation}
(див., напр., \cite[означення~12.1]{Va}). У подальшому
\begin{equation}\label{eq47***}
h(E)=\sup\limits_{x, y\in E}h(x, y)
\end{equation}
 -- хордальний діаметр множини
$E\subset \overline{{\Bbb R}^n}.$ Як звично, сім'я $\frak{F}$
відображень $f:D\rightarrow \overline{{\Bbb C}}$ буде називатися
{\it нормальною,} якщо з кожної послідовності $f_n\in \frak{F},$
$n=1,2,\ldots , $ можна виділити підпослідовність $f_{n_k},$
$k=1,2,\ldots ,$ яка збігається локально рівномірно до деякого
відображення $f:D\rightarrow \overline{{\Bbb C}}$ в метриці $h.$
Якщо додатково $f\in \frak{F},$ сім'я $\frak{F}$ називається {\it
компактною}.

\medskip
Множина $A\subset {\Bbb D}$ називається {\it інваріантно опуклою,}
якщо множина $g(A)$ є опуклою для будь-якого дробово-лінійного
автоморфізму $g$ одиничного круга.

\medskip
Нехай $D$ -- область в ${\Bbb R}^n.$ Будемо говорити, що функція
${\varphi}:D\rightarrow{\Bbb R},$ що є локально інтегровною в
деякому околі точки $x_0\in D,$ має {\it скінченне середнє
коливання} в точці $x_0$ (пишемо: $\varphi\in FMO(x_0)$), якщо
\begin{equation}\label{eq17:}
{\limsup\limits_{\varepsilon\rightarrow
0}}\frac{1}{\Omega_n\varepsilon^n}\int\limits_{B(
x_0,\,\varepsilon)}
|{\varphi}(x)-\overline{{\varphi}}_{\varepsilon}|\ dm(x)\, <\,
\infty\,,
\end{equation}
де $\Omega_n$ -- об'єм одиничної кулі в ${\Bbb R}^n,$
$\overline{{\varphi}}_{\varepsilon}=\frac{1}{\Omega_n\varepsilon^n}\int\limits_{B(
x_0,\,\varepsilon)} {\varphi}(x)\ dm(x)$ (див., напр.,
\cite[розд.~2]{RSY$_2$}).
Зауважимо, що коли виконується умова $(\ref{eq17:})$ можлива
ситуація, коли $\overline{{\varphi}_{\varepsilon}}\rightarrow\infty$
при $\varepsilon\rightarrow 0.$
Також будемо говорити, що  ${\varphi}:D\rightarrow{\Bbb R}$ --
функція скінченного середнього коливання \ в \ області\  D, пишемо
${\varphi}\in FMO(D),$ якщо ${\varphi}$ має скінченне середнє
коливання в кожній точці $x_0\in D.$

\medskip Нехай $M(z)\subset {\Bbb D},$ $z\in {\Bbb C}$ -- деяка система множин (тобто,
при кожному $z_0\in {\Bbb C}$ символ $M(z_0)$ позначає деяку множину
в ${\Bbb D}$). Позначимо через $\frak{M}_M$ множину всіх комплексних
вимірних функцій $\mu:{\Bbb C}\rightarrow {\Bbb D},$ таких що
$\mu(z)\in M(z)$ при майже всіх $z\in {\Bbb C}.$ Нехай $K$ --
компакт в ${\Bbb C},$ $M(z)$ -- система множин в ${\Bbb D}.$
Позначимо через $\frak{F}_M(K)$ клас усіх регулярних розв'язків
$f:{\Bbb C}\rightarrow{\Bbb C}$ рівняння~(\ref{eq2}) з комплексними
коефіцієнтами $\mu$, рівними нулю зовні $K$ такими, що
\begin{equation}\label{eq1C}
f(z)=z+o(1)\quad {\text при}\quad z\rightarrow\infty\,,
\end{equation}
при цьому $\mu\in \frak{M}_M.$ Для множини $M(z)$ покладемо
\begin{equation}\label{eq1K} Q_M(z)=\frac{1+q_M(z)}{1-q_M(z)}\,,\quad
q_M(z)=\sup\limits_{\nu\in M(z)}|\nu|\,.
\end{equation}
Одним з основних результатів статті є наступне твердження.
\medskip
\begin{theorem}\label{th1}
{\sl\, Нехай $M(z),$ $z\in {\Bbb C}$ -- сім'я інваріантно опуклих
компактних множин, і нехай функція $Q_M$ є інтегровною на $K$ і
задовольняє принаймні одну з умов: або $Q_M\in FMO({\Bbb C}),$ або
для кожного $z_0\in {\Bbb C}$ існує $\delta_0=\delta(z_0)>0$ таке,
що
\begin{equation}\label{eq2D}
\int\limits_0^{\delta_0}\frac{dt}{tq_{M_{z_0}}(t)}=\infty\,,
\end{equation}
де
$q_{M_{z_0}}(t)=\frac{1}{2\pi}\int\limits_{0}^{2\pi}Q_M(z_0+e^{it})\,dt.$
Тоді сім'я відображень $\frak{F}_M(K)$ є компактною в ${\Bbb C}.$}
\end{theorem}

\medskip
Перейдемо тепер до розгляду питання про компактність класів
розв'язків задачі Діріхле для рівняння Бельтрамі. Розглянемо
наступну задачу Діріхле:
\begin{equation}\label{eq2C}
f_{\overline{z}}=\mu(z)\cdot f_z\,,
\end{equation}
\begin{equation}\label{eq1A}
\lim\limits_{\zeta\rightarrow z}{\rm
Re\,}f(\zeta)=\varphi(z)\qquad\forall\,\, z\in \partial D\,,
\end{equation}
де $\varphi:\partial D\rightarrow {\Bbb R}$ -- наперед задана
неперервна функція. Надалі вважаємо, що $D$ -- деяка однозв'язна
жорданова область у ${\Bbb C}.$ Розв'язок
задачі~(\ref{eq2C})--(\ref{eq1A}) будемо вважати {\it регулярним,}
якщо виконано одно з двох: або $f(z)=const$ в $D,$ або $f$ --
відкрите дискретне відображення класу $W_{\rm loc}^{1, 1}(D),$ таке
що $J(z, f)\ne 0$ при майже всіх $z\in D.$

\medskip
Зафіксуємо точку $z_0\in D$ і функцію $\varphi.$ Нехай $M(z)\subset
{\Bbb D},$ $z\in D$ -- деяка система множин. Позначимо через
$\frak{M}_M$ множину всіх комплексних вимірних функцій
$\mu:D\rightarrow {\Bbb D},$ таких що $\mu(z)\in M(z)$ при майже
всіх $z\in D.$ Нехай $\frak{F}_{\varphi, M, z_0}(D)$ позначає клас
усіх регулярних розв'язків $f:D\rightarrow{\Bbb C}$ задачі
Діріхле~(\ref{eq2C})--(\ref{eq1A}), які задовольняють умову ${\rm
Im}\,f(z_0)=0$ таких, що $\mu\in \frak{M}_M.$ Як і раніше, визначимо
функцію $Q_M(z)$ співвідношенням~(\ref{eq1K}), причому вважатимемо
$Q_M(z)\equiv 1$ при $z\in {\Bbb C}\setminus D.$ Наступне твердження
узагальнює~\cite[теорема~2]{Dyb} на випадок довільних однозв'язних
жорданових областей.

\medskip
\begin{theorem}\label{th2A}
{\sl Нехай $D$ -- деяка однозв'язна жорданова область у ${\Bbb C},$
і нехай функція $\varphi(z)$ у~(\ref{eq1A}) неперервна. Припустимо,
що $M(z),$ $z\in D,$ -- сім'я інваріантно опуклих компактних множин,
і нехай функція $Q_M$ є інтегровною в $D$ і задовольняє принаймні
одну з умов: або $Q_M\in FMO(\overline{D}),$ або для кожного $x_0\in
\overline{D}$ існує $\delta_0=\delta(x_0)>0$ таке, що
\begin{equation}\label{eq2G}
\int\limits_0^{\delta_0}\frac{dt}{tq_{M_{x_0}}(t)}=\infty\,,
\end{equation}
де
$q_{M_{x_0}}(t)=\frac{1}{2\pi}\int\limits_0^{2\pi}Q_M(x_0+e^{it})\,dt.$
Тоді сім'я відображень $\frak{F}_{\varphi, M, z_0}(D)$ є компактною
в $D.$}
\end{theorem}

\medskip
{\bf 2. Про збіжність гомеоморфізмів з модульними умовами.} Нехай
$D$ -- область в ${\Bbb R}^n,$ $n\geqslant 2,$ $M(\Gamma)$ --
конформний модуль сім'ї кривих $\Gamma$ в ${\Bbb R}^n$ (див., напр.,
\cite[гл.~6]{Va}). Покладемо
$$S(x_0, r)=\{x\in {\Bbb R}^n: |x-x_0|=r\}\,,B(x_0, r)=\{x\in {\Bbb R}^n:
|x-x_0|<r\}\,,$$
$${\Bbb B}^n:=B(0, 1)\,,\quad {\Bbb S}^{n-1}:=S(0, 1)\,,\quad \Omega_n=m({\Bbb B}^n)\,,\quad
\omega_{n-1}=\mathcal{H}^{n-1}({\Bbb S}^{n-1})\,,$$
де $\mathcal{H}^{n-1}$ позначає $(n-1)$-вимірну міру Хаусдорфа в
${\Bbb R}^n.$ Якщо $n=2,$ покладемо ${\Bbb D}:=B(0, 1).$ Нехай, крім
того,
$$A=A(x_0, r_1, r_2)=\{x\in {\Bbb R}^n: r_1<|x-x_0|<r_2\}\,.$$
Для заданих множин $E,$ $F\subset\overline{{\Bbb R}^n}$ і області
$D\subset {\Bbb R}^n$ позначимо через $\Gamma(E,F,D)$ сім'ю всіх
кривих $\gamma:[a,b]\rightarrow \overline{{\Bbb R}^n}$ таких, що
$\gamma(a)\in E,\gamma(b)\in\,F$ і $\gamma(t)\in D$ при $t \in [a,
b].$ Відображення $f:D\rightarrow \overline{{\Bbb R}^n}$ називається
{\it кільцевим $Q$-ві\-доб\-ра\-жен\-ням у точці
$x_0\,\in\,\overline{D},$} якщо співвідношення
\begin{equation}\label{eq3*!gl0} M(f(\Gamma(S(x_0, r_1),\,S(x_0, r_2),\,A\cap D)))
\leqslant \int\limits_{A\cap D} Q(x)\cdot \eta^n(|x-x_0|)\, dm(x)
\end{equation}
виконано для будь-якого кільця $A=A(x_0, r_1,r_2),$ $0<r_1<r_2<
r_0:=\sup\limits_{x\in D}|x-x_0|,$ і кожної вимірної за Лебегом
функції $\eta:(r_1,r_2)\rightarrow [0,\infty ]$ такої, що
\begin{equation}\label{eq*3gl0}\int\limits_{r_1}^{r_2}\eta(r)\,dr \geqslant 1\,.
\end{equation}

\medskip
Наступна важлива лема випливає з~\cite[теореми~4.1 і 4.2]{RSS}.

\medskip
\begin{lemma}\label{lem1}
{\sl\, Нехай $D$ -- область в ${\Bbb R}^n,$ $n\geqslant 2,$
$Q:D\rightarrow [1, \infty]$ -- вимірна за Лебегом функція. Нехай,
крім того, $f_k,$ $k=1,2,\ldots $ -- послідовність гомеоморфізмів
області $D$ в ${\Bbb R}^n,$ які задовольняють
умови~(\ref{eq3*!gl0})--(\ref{eq*3gl0}) в кожній точці $x_0$ області
$D,$ що збігається локально рівномірно в $D$ до деякого відображення
$f:D\rightarrow \overline{{\Bbb R}^n}$ по хордальній метриці $h.$
Припустимо, що функція $Q$ задовольняє принаймні одну з двух умов:
або $Q\in FMO(D),$ або для кожного $x_0\in D$ існує
$\delta_0=\delta(x_0)>0$ таке, що
\begin{equation}\label{eq2H}
\int\limits_0^{\delta_0}\frac{dt}{tq^{\frac{1}{n-1}}_{x_0}(t)}=\infty\,,
\end{equation}
де $q_{x_0}(t)=\frac{1}{\omega_{n-1}r^{n-1}}\int\limits_{S(x_0,
t)}Q(x)\,d\mathcal{H}^{n-1},$ а $\mathcal{H}^{n-1}$ позначає
$(n-1)$-вимірну міру Хаусдорфа. Тоді відображення $f$ є або
гомеоморфізмом $f:D\rightarrow {\Bbb R}^n,$ або сталою
$c\in\overline{{\Bbb R}^n}.$ }
\end{lemma}

\medskip
Згідно з~\cite{GM}, область $D$ в ${\Bbb R}^n$ називається {\it
областю квазіекстремальної довжини}, скор. {\it $QED$-об\-лас\-тю},
якщо знайдеться число $A\geqslant 1$ таке, що для всіх континуумів
$E$ і $F$ у $D$ виконується нерівність
\begin{equation}\label{eq4***}
M(\Gamma(E, F, {\Bbb R}^n))\leqslant A\cdot M(\Gamma(E, F, D))\,.
\end{equation}

\medskip
Наступне твердження встановлено в~\cite[лема~3.2]{Sev$_4$}.

\medskip
\begin{lemma}\label{lem2}
{\sl\, Нехай $D,$ $D^{\,\prime}$ -- області в ${\Bbb R}^n,$
$n\geqslant 2,$ $b_0\in D,$ $b_0^{\,\prime}\in D^{\,\prime},$ і
нехай $f_k,$ $k=1,2,\ldots,$ -- сім'я гомеоморфізмів області $D$ на
$D^{\,\prime}$ з умовою нормування $f_k(b_0)=b^{\,\prime}_0,$
$k=1,2,\ldots .$ Припустимо, що кожне відображення $f_k,$
$k=1,2,\ldots$ задовольняє співвідношення~(\ref{eq3*!gl0}) в
довільній точці $x_0\in \overline{D}$ і деякою вимірною функцією
$Q:D\rightarrow [1, \infty]$ такою, що виконано принаймні одну з
двох умов: або $Q\in FMO(\overline{D}),$ або для кожного $x_0\in
\overline{D}$ існує $\delta_0=\delta(x_0)>0$ таке, що
виконується~(\ref{eq2H}). Нехай область $D$ є локально зв'язною на
своїй межі, а $D^{\,\prime}$ є $QED$-областю, яка містить не менше
одної скінченної межової точки. Тоді кожне $f_k,$ $k=1,2,\ldots,$
продовжується до неперервного відображення
$\overline{f}_k:\overline{D}\rightarrow\overline{D^{\,\prime}}$ і,
крім того, сім'я продовжених відображень $\overline{f}_k,$
$k=1,2,\ldots,$ є одностайно неперервною в $\overline{D}.$
 }
\end{lemma}

\medskip
{\bf 3. Компактність розв'язків рівнянь Бельтрамі з гідродинамічним
нормуванням.} {\it Доведення теореми~\ref{th1}}. {\bf I.} Передусім
доведемо, що сім'я $\frak{F}_M(K)$ є одностайно неперервною.
Зафіксуємо $f\in\frak{F}_M(K),$ довільний компакт $C\subset {\Bbb
C}$ і покладемо $\widetilde{f}=\frac{1}{f(1/z)}.$ (Якщо
$f(1/w_0)=0,$ то ми вважаємо, що $\widetilde{f}(w_0):=\infty$).
Оскільки $f(z)=z+o(1)$ при $z\rightarrow\infty,$ то
$\lim\limits_{z\rightarrow\infty}f(z)=\infty.$ Тоді покладаючи
$\widetilde{f}(0):=0,$ отримаємо гомеоморфізм $\widetilde{f},$
визначений в деякому околі нуля. У подальшому покладемо
$f(\infty):=\infty.$ Оскільки $f(z)=z+o(1)$ при
$z\rightarrow\infty,$ існує окіл $U$ початку координат і функція
$\varepsilon: U\rightarrow {\Bbb C}$ такі, що
$f(1/z)=1/z+\varepsilon(1/z),$ де $z\in U$ і
$\varepsilon(1/z)\rightarrow 0$ при $z\rightarrow 0.$ Отже, при
достатньо малому $\Delta z\in {\Bbb C}$ ми будемо мати, що
$$\frac{\widetilde{f}(\Delta z)-\widetilde{f}(0)}{\Delta z}=
\frac{1}{\Delta z}\cdot \frac{1}{1/(\Delta z)+\varepsilon(1/\Delta
z)}=\frac{1}{1+(\Delta z)\cdot \varepsilon(1/\Delta z)}\rightarrow 1
$$
при $z\rightarrow 0.$ Сказане доводить, що існує
$\widetilde{f}^{\,\prime}(0),$ при цьому,
$\widetilde{f}^{\,\prime}(0)=1.$ Оскільки зовні $K$ функція $\mu$
дорівнює нулю, відображення $f$ є конформним в деякому околі
$V:={\Bbb C}\setminus B(0, 1/r_0)$ точки нескінченність, причому,
число $1/r_0$ залежить тільки від $K$ і $K\subset B(0, 1/r_0).$ Без
обмеження загальності можна вважати, що також і компакт $C$
задовольняє умову $C\subset B(0, 1/r_0).$   В такому випадку,
відображення $\widetilde{f}=\frac{1}{f(1/z)}$ є конформним у кулі
$B(0, r_0),$ при цьому відображення $F(z):=\frac{1}{r_0}\cdot
\widetilde{f}(r_0z)$ є гомеоморфізмом одиничного круга в ${\Bbb C}$
і задовольняє умови $F(0)=0,$ $F^{\,\prime}(1)=1.$ За теоремою Кебе
про 1/4 (див. напр., \cite[теорема, розд.~1.3 гл.~1]{GR$_1$})
$F({\Bbb D})\supset B(0, 1/4);$ тоді
\begin{equation}\label{eq4}
\widetilde{f}(B(0, r_0))\supset B(0, r_0/4)\,.
\end{equation}
Зі співвідношення~(\ref{eq4}) випливає, що
\begin{equation}\label{eq5}
(1/f)(\overline{\Bbb C}\setminus \overline{B(0, 1/r_0)})\supset B(0,
r_0/4)\,.
\end{equation}
З урахуванням~(\ref{eq5}) покажемо, що
\begin{equation}\label{eq6}
f(\overline{\Bbb C}\setminus \overline{B(0, 1/r_0)})\supset
\overline{{\Bbb C}}\setminus \overline{B(0, 4/r_0)}\,.
\end{equation}
Дійсно, нехай $y\in \overline{{\Bbb C}}\setminus \overline{B(0,
4/r_0)},$ тоді $\frac{1}{y}\in B(0, r_0/4).$ Зі
співвідношення~(\ref{eq5}) $\frac{1}{y}=(1/f)(x),$ $x\in
\overline{\Bbb C}\setminus \overline{B(0, 1/r_0)}.$ Тоді $y=f(x),$
$x\in \overline{\Bbb C}\setminus \overline{B(0, 1/r_0)},$ що і
доводить~(\ref{eq6}).

\medskip
Оскільки $f$ -- гомеоморфізм у ${\Bbb C},$ зі
співвідношення~(\ref{eq6}) випливає, що $f(B(0, 1/r_0))\subset B(0,
4/r_0).$ Покладемо $\Delta:=h(\overline{\Bbb C}\setminus B(0,
4/r_0)),$ де $h(\overline{\Bbb C}\setminus B(0, 4/r_0))$ --
хордальний діаметр множини $\overline{\Bbb C}\setminus B(0, 4/r_0).$
За~\cite[теорема~3.1]{LSS} кожне відображення $f$ є так званим
кільцевим $Q$-гомеоморфізмом в ${\Bbb C}$ при $Q=K_{\mu}(z),$ де
$\mu$ визначається зі співвідношення~(\ref{eq2}), а $K_{\mu}$
визначено в~(\ref{eq1}). Зауважимо, що $Q\leqslant Q_M(z)$ майже
скрізь. В таком випадку, сім'я відображень $\frak{F}_M(K)$ є
одностайно неперервною в~$B(0, 1/r_0)$ за~\cite[теореми~7.5,
7.6]{MRSY}. Нехай тепер $f_n\in \frak{F}_M(K),$ $n=1,2,\ldots .$ За
теоремою Арцела-Асколі (див., напр., \cite[теорема~20.4]{Va}) існує
підпослідовність $f_{n_k}(z)$ послідовності $f_n,$ $k=1,2,\ldots ,$
а також неперервне відображення $f:B(0,
1/r_0)\rightarrow\overline{\Bbb C}$ такі, що $f_{n_k}$ локально
рівномірно збігається до відображення $f$ у $B(0, 1/r_0)$ при
$k\rightarrow\infty.$ Зокрема, оскільки компакт $C$ належить $B(0,
1/r_0),$ послідовність $f_{n_k}$ збігається до $f$ рівномірно на
$C.$ Оскільки компакт $C$ був обраний довільним, ми встановили, що
сім'я відображень $f_{n_k}$ збігається до відображення $f$ локально
рівномірно.

\medskip
{\bf II.} Для завершення доведення теореми~\ref{th1} залишилось
встановити, що $f\in \frak{F}_M(K).$ Передусім доведемо, що граничне
відображення $f$ задовольняє умову $f(z)=z+o(1)$ при
$z\rightarrow\infty.$ Зауважимо, що сім'я відображень
$F_{n_k}(z):=\frac{1}{r_0}\cdot\frac{1}{f_{n_k}(\frac{1}{r_0z})}$ є
компактною в одиничному крузі (див.~\cite[теорема~1.2
гл.~I]{GR$_1$}). Без обмеження загальності можна вважати, що сама
послідовність $F_{n_k}$ є локально рівномірно збіжною в ${\Bbb D}.$
Тоді функція $F(z)=\frac{1}{r_0}\cdot\frac{1}{f(\frac{1}{r_0z})}$
знову належить до класу $S,$ що складається з конформних відображень
$F$ одиничного круга, які задовольняють умови $F(0)=0,$
$F^{\,\prime}(0)=1.$ Розпишемо розклад в ряд Тейлора функцій
$F_{n_k}$ і $F$ в околі нуля. Будемо мати:
\begin{equation}\label{eq1E}F_{n_k}(z)=z+c_kz^2+z^2\cdot \varepsilon_k(z)\,,\quad
k=1,2,\ldots\,,
\end{equation}
\begin{equation}\label{eq2K}
F(z)=z+c_0z^2+z^2\cdot \varepsilon_0(z)\,,
\end{equation}
де $\varepsilon_k(z)$ і $\varepsilon_0(z)$ прямують до нуля при
$z\rightarrow 0.$ Зі співвідношень~(\ref{eq1E}) і (\ref{eq2K})
випливає, що
\begin{equation}\label{eq3A}
f_{n_k}(t)=\frac{r_0t^2}{r_0t+c_k+\varepsilon_k\left(\frac{1}{r_0t}\right)}\,,\quad
f_{n_k}(t)-t=-\frac{c_k+\varepsilon_k\left(\frac{1}{r_0t}\right)}
{r_0+\frac{c_k}{t}+\frac{\varepsilon_k\left(\frac{1}{r_0t}\right)}{t}}\,,
\quad k=1,2,\ldots\,,
\end{equation}
\begin{equation}\label{eq4B}
f(t)=\frac{r_0t^2}{r_0t+c_0+\varepsilon_0\left(\frac{1}{r_0t}\right)}\,,\quad
f(t)-t=-\frac{c_0+\varepsilon_0\left(\frac{1}{r_0t}\right)}
{r_0+\frac{c_0}{t}+\frac{\varepsilon_0\left(\frac{1}{r_0t}\right)}{t}}\,.
\end{equation}
Зокрема, з другого співвідношення у~(\ref{eq3A}) переходом до
границі при $t\rightarrow\infty$ випливає, що
$f_{n_k}(t)-t\rightarrow -\frac{c_k}{r_0}.$ Оскільки за умовою
$f_{n_k}(t)-t\rightarrow 0$ при $t\rightarrow \infty,$ маємо:
$c_k=0.$ За теоремою Вейерштрасса про збіжність коефіцієнтів ряду
Тейлора (див., напр., \cite[Теорема~1.1.I]{Gol}) маємо:
$c_k=0\rightarrow c_0$ при $k\rightarrow\infty.$ Отже, $c_0=0$ в
(\ref{eq4B}), тобто, відображення $f$ також має гідродинамічне
нормування: $f(z)=z+o(1)$ при $z\rightarrow\infty.$

\medskip
Тепер покажемо, що $f$ -- гомеоморфізм комплексної площини.
Покладемо $\mu_k:=\mu_{f_{n_k}}.$ За~\cite[теорема~3.1]{LSS} кожне
відображення $f_{n_k}$ є кільцевим $Q$-гомеоморфізмом в кожній точці
$z_0\in {\Bbb C}$ при $Q=K_{\mu_k}(z).$ З огляду на лему~\ref{lem1}
має місце наступна альтернатива: або $f$ -- гомеоморфізм з $D$ у
${\Bbb C},$ або $f$ -- стала в $\overline{\Bbb C}.$ За доведеним на
кроці~{\bf I} $f$ є гомеоморфізмом в деякому околі нескінченності,
отже, $f$ -- гомеоморфізм всієї комплексної площини, який приймає
тільки скінченні комплексні значення.

\medskip
Тоді за~\cite[теорема~3.1 і зауваження~3.1]{RSY$_3$} $f\in W_{\rm
loc}^{1, 1}({\Bbb C}).$ Покажемо, що $f$ є регулярним розв'язком
рівняння~(\ref{eq2}) з деяким $\mu$, тобто, $J(z, f)\ne 0$ при майже
всіх $z\in {\Bbb C}.$ Оскільки за умовою відображення $f_{n_k}$ є
регулярними, крім того, $K_{\mu_k}(z)\leqslant Q_M(z)$ майже скрізь,
крім того, $Q_M$ є інтегровною в $K,$ маємо: $J(f_{n_k}, z)\ne 0$
при майже всіх $z\in {\Bbb C},$ $g_{n_k}:=f^{\,-1}_{n_k}\in W_{\rm
loc}^{1, 2},$ причому, $(g_{n_k})_w=0$ у майже всіх точках $w,$ де
$J(g_{n_k}, w)=0$ (див. \cite[теорема~1.3]{HK$_1$}). Отже, $g_{n_k}$
мають $N$-властивість Лузіна (див., напр., \cite[теорема~1]{Pon} або
\cite[наслідок~B]{MM}). Позначимо $\partial g=g_w$ і
$\overline{\partial}g=g_{\overline{w}}.$ Нехай $C_0$ -- довільний
компакт у $f(D).$ Оскільки $f_{n_k}$ збігаються до $f$ локально
рівномірно в ${\Bbb C}$ і $f$ -- гомеоморфізм, то $g_{n_k}$ також
збігаються локально рівномірно в $f(D)$ до відображення
$g:=f^{\,-1}$ (див., напр., \cite[лема~3.1]{RSS}). Зокрема, звідси
випливає, що відображення $g_{n_k}$ визначені на $C_0$ при всіх
достатньо великих $k\geqslant k_0=k_0(C_0).$ Оскільки $g_{n_k}$
збігаються рівномірно на $C_0$ до $g$ при $k\rightarrow\infty,$
маємо:
\begin{equation}\label{eq7B}|g_{n_k}(w)|\leqslant |g(w)|+1
\end{equation}
при всіх $w\in C_0$ і достатньо великих $k\in {\Bbb N}.$ Покладемо
$A:=\sup\limits_{w\in C_0}|g(w)|.$ Зауважимо, що за теоремою Кантора
про обмеженість неперервного відображення на компакті $C_0$
виконується нерівність $A<\infty.$ Отже, з огляду на~(\ref{eq7B}),
\begin{equation}\label{eq8B}
g_{n_k}(C_0)\subset B(0, A+1)\,.
\end{equation}
Розглянемо спочатку випадок, коли $J(g_{n_k}, w)\ne 0$ при майже
всіх $w\in C_0.$ Будемо мати:
$$|\partial g_{n_k}(w)|^2=\left(|\partial
g_{n_k}(w)|^{\,2}-|\overline{\partial}g_{n_k}(w)|^{\,2}\right)\cdot
\frac{|\partial g_{n_k}(w)|^{\,2}}{\left(|\partial
g_{n_k}(w)|^{\,2}-|\overline{\partial}g_{n_k}(w)|^{\,2}\right)}=$$$$=J(w,
g_{n_k})\cdot\frac{1}{1-\left|\frac{\overline{\partial}g_{n_k}(w)}{\partial
g_{n_k}(w)}\right|^2}\,,$$
зокрема,
\begin{equation}\label{eq3}
|\partial g_{n_k}(w)|^2=J(w,
g_{n_k})\cdot\frac{1}{1-\left|\frac{\overline{\partial}g_{n_k}(w)}{\partial
g_{n_k}(w)}\right|^2}\,.
\end{equation}
Оскільки $g_{n_k}$ є гомеоморфізмами класу $W_{\rm loc}^{1, 1}(D),$
які мають $N$-властивість Лузіна, для них має місце заміна змінних в
інтегралі (див., напр.,~\cite[теорема~3.2.5]{Fe}). Зауважимо, що
$\mu_{g_{n_k}}(w)=-\mu_k(g_{n_k}(w))=-\mu_k(f^{\,-1}_{n_k}(w)),$
див. напр.,~\cite[(4).C.I]{A}. В такому випадку, з огляду
на~(\ref{eq8B}) і (\ref{eq3}), ми будемо мати, що
$$\int\limits_{C_0}|\partial g_{n_k}(w)|^2\,dm(w)=
\int\limits_{C_0}\left(|\partial
g_{n_k}(w)|^{\,2}-|\overline{\partial}g_{n_k}(w)|^{\,2}\right)\cdot
\frac{|\partial g_{n_k}(w)|^{\,2} dm(w)}{\left(|\partial
g_{n_k}(w)|^{\,2}-|\overline{\partial}g_{n_k}(w)|^{\,2}\right)}=$$
\begin{equation}\label{eq4F}
=\quad \int\limits_{C_0} J(w,
g_{n_k})\cdot\frac{1}{1-\left|\frac{\overline{\partial}g_{n_k}(w)}{\partial
g_{n_k}(w)}\right|^2}\,\,dm(w)= \quad\int\limits_{g_{n_k}(C_0)}
\frac{dm(z)}{1-|\mu_k(z)|^2}\quad\leqslant
\end{equation}
$$\leqslant \quad\int\limits_{B(0, A+1)} Q^{\,\prime}(z)\, dm(z)<\infty\,.$$
З~(\ref{eq4F}) випливає, що
\begin{equation}\label{eq9B}
\int\limits_{C_0}|\partial g_{n_k}(w)|^2\,dm(w)\leqslant
\int\limits_{B(0, A+1)} Q^{\,\prime}(z)\, dm(z)<\infty\,.
\end{equation}
Зауважимо, що співвідношення~(\ref{eq9B}) виконується також і у
випадку, коли $J(g_{n_k}, w)=0$ на множині додатної міри в $C_0,$
бо, як вже було зауважено вище, $(g_{n_k})_w=0$ у майже всіх точках
$w,$ де $J(g_{n_k}, w)=0$ (див. \cite[теорема~1.3]{HK$_1$}).
З~(\ref{eq4F}) випливає, що $g\in W_{\rm loc}^{1, 2}$
(див.~\cite[лема~III.3.5]{Re}). Тоді $g$ має $N$-властивість за
теоремою Малого-Мартіо, див., напр., \cite[наслідок~B]{MM}. У свою
чергу, $f$ має майже всюди невироджений якобіан за теоремою
Пономарьова (див. \cite[теорема~1]{Pon}).

\medskip
Отже, $f$ є регулярним розв'язком рівняння~(\ref{eq2}) при деякій
функції $\mu:{\Bbb C}\rightarrow {\Bbb D}.$ За теоремою
Герінга-Лехто відображення $f$ є майже всюди диференційовним
(див.~\cite[теорема~III.3.1]{LV}). Тому за теоремою~16.1
у~\cite{RSS} $\mu(z)=0$ при всіх $z\in {\Bbb C}\setminus K.$
Нарешті, $\mu\in\frak{M}_M$ за \cite[лема~1]{L$_2$}. Отже, $f\in
\frak{F}_M(K).$~$\Box$

\medskip
{\bf 4. Одностайна неперервність сімей відображень з оберненою
нерівністю Полецького.} В цьому розділі ми маємо справу з
відображеннями $f:D\rightarrow{\Bbb R}^n$ області $D\subset{\Bbb
R}^n,$ $n\geqslant 2.$ Нижче ми наводимо деякі результати
з~\cite{SSD}, що є необхідним для доведення ключових теорем
наступного розділу.

\medskip
Нехай, як і раніше, $M(\Gamma)$ -- конформний модуль сім'ї кривих
$\Gamma$ в ${\Bbb R}^n$ (див., напр., \cite[гл.~6]{Va}). Якщо
$f:D\rightarrow {\Bbb R}^n$ -- задане відображення, $y_0\in f(D)$ і
$0<r_1<r_2<d_0=\sup\limits_{y\in f(D)}|y-y_0|,$ то через
$\Gamma_f(y_0, r_1, r_2)$ позначимо сім'ю всіх кривих $\gamma$ в
області $D$ таких, що $f(\gamma)\in \Gamma(S(y_0, r_1), S(y_0, r_2),
A(y_0,r_1,r_2)).$ Нехай $Q:{\Bbb R}^n\rightarrow [0, \infty]$ --
вимірна за Лебегом функція.  Будемо говорити, що {\it $f$
задовольняє обернену нерівність Полецького} в точці $y_0\in f(D),$
якщо співвідношення
\begin{equation}\label{eq2*A}
M(\Gamma_f(y_0, r_1, r_2))\leqslant \int\limits_{A(y_0,r_1,r_2)\cap
f(D)} Q(y)\cdot \eta^n (|y-y_0|)\, dm(y)
\end{equation}
виконується для довільної вимірної за Лебегом функції $\eta:
(r_1,r_2)\rightarrow [0,\infty ]$ такої, що
\begin{equation}\label{eqA2}
\int\limits_{r_1}^{r_2}\eta(r)\, dr\geqslant 1\,.
\end{equation}
Нагадаємо, що відображення $f:D\rightarrow D^{\,\prime}$ називається
{\it замкненим,} якщо будь-яку замкнену множину $A\subset D$
відображення $f$ переводить у замкнену множину $f(A)\subset
D^{\,\prime}$ (де замкненість розуміється відносно областей $D$ і
$D^{\,\prime},$ відповідно). Можна показати, що будь-які
гомеоморфізми між областями $D$ і $D^{\,\prime}$ є замкненими
відображеннями.

\medskip Для областей $D, D^{\,\prime}\subset {\Bbb R}^n,$ точок
$a\in D,$ $b\in D^{\,\prime}$ і вимірної за Лебегом функції $Q\colon
D\rightarrow [0, \infty],$ що дорівнює нулю зовні області $D,$
позначимо через ${\frak S}_{a, b, Q}(D, D^{\,\prime})$ сім'ю всіх
відкритих, дискретних і замкнених відображень $f$ області $D$ на
$D^{\,\prime},$ що задовольняють умову~(\ref{eq2*A}) для кожного
$y_0\in f(D),$ таких що $f(a)=b.$ Наступне твердження доведено
в~\cite[теорема~7.1]{SSD}.

\medskip
\begin{lemma}\label{th4}
{\sl Припустимо, що область $D$ має слабо плоску межу, жодна із
зв'язних компонент якої не вироджена. Нехай $Q\in L^1(D).$ Якщо
область $D^{\,\prime}$ локально зв'язна на своїй межі, то кожне
відображення $f\in {\frak S}_{a, b, Q}(D, D^{\,\prime})$ має
неперервне продовження $\overline{f}:\overline{D}\rightarrow
\overline{D^{\,\prime}},$ причому,
$\overline{f}(\overline{D})=\overline{D^{\,\prime}}$ і, крім того,
сім'я ${\frak S}_{a, b, Q}(\overline{D}, \overline{D^{\,\prime}})$
усіх продовжених відображень $\overline{f}:\overline{D}\rightarrow
\overline{D^{\,\prime}}$ є одностайно неперервною в $\overline{D}.$
}
\end{lemma}

\medskip
{\bf 5. Компактність розв'язків задачі Діріхле.} {\it Доведення
теореми~\ref{th2A}}. {\textbf I.} Нехай $f_m,$ $m=1,2,\ldots $ --
довільна послідовність сім'ї $\frak{F}_{\varphi, M, z_0}(D).$ Згідно
теореми Стойлова про факторизацію (див., напр., \cite[п.~5 (III),
гл.~V]{St}), для відображення $f_m$ справедливе зображення
\begin{equation}\label{eq2E}
f_m=\varphi_m\circ g_m\,,
\end{equation}
де $g_m$ -- деякий гомеоморфізм, а $\varphi_m$ -- аналітична
функція. За лемою~1 в~\cite{Sev$_3$} відображення $g_m$ належить
класу Соболєва $W_{\rm loc}^{1, 1}(D)$ і має скінченне спотворення.
Більше того, згідно~\cite[(1), п.~C, гл.~I]{A} для майже всіх $z\in
D$ отримаємо:
\begin{equation}\label{eq1B}
{f_m}_z={\varphi_m}_z(g_m(z)){g_m}_z,\qquad
{f_m}_{{\overline{z}}}={\varphi_m}_z(g_m(z)){g_m}_{\overline{z}}\,.
\end{equation}
Отже, за співвідношенням~(\ref{eq1B}), $J(z, g_m)\ne 0$ для майже
всіх $z\in D,$ крім того, $K_{\mu_{f_m}}(z)=K_{\mu_{g_m}}(z).$

\medskip
\textbf{II.} Доведемо, що межа області $g_m(D)$ містить не менше
двох точок. Припустимо супротивне. Тоді або $g_m(D)={\Bbb C},$ або
$g_m(D)={\Bbb C}\setminus \{a\},$ де $a\in {\Bbb C}.$ Нехай спочатку
$g_m(D)={\Bbb C}.$ Тоді за теоремою Пікара $\varphi_m(g_m(D))$ є
всією площиною, за виключенням, можливо, однієї точки $\omega_0\in
{\Bbb C}.$ З іншого боку, при кожному $m=1,2,\ldots$ функція
$u_m(z):={\rm Re}\,f_m(z)={\rm Re}\,(\varphi_m(g_m(z)))$ неперервна
на компакті $\overline{D}$ за умовою~(\ref{eq1A}) і з огляду на
неперервність функції $\varphi.$ Отже, існує $C_m>0$ таке, що $|{\rm
Re}\,f_m(z)|\leqslant C_m,$ але це суперечить тому, що
$\varphi_m(g_m(D))$ містить всі точки комплексної площини крім,
можливо, однієї. Випадок $g_m(D)={\Bbb C}\setminus \{a\},$ $a\in
{\Bbb C},$ є неможливим, оскільки $g_m(D)$ має бути однозв'язною
областю в ${\Bbb C}$ як образ однозв'язної області $D$ при
гомеоморфізмі $g_m,$ $m=1,2,\ldots \,.$

\medskip
Отже,  межа області $g_m(D)$ містить не менше двох точок. Тоді за
теоремою Рімана про відображення можна перетворити область $g_m(D)$
на одиничний круг ${\Bbb D}$ за допомогою конформного відображення
$\psi_m.$ Нехай $z_0\in D$ -- точка з умови теореми. За допомогою
допоміжного конформного відображення
$$\widetilde{\psi_m}(z)=\frac{z-(\psi_m\circ
g_m)(z_0)}{1-z\overline{(\psi_m\circ g_m)(z_0)}}$$ одиничного круга
на себе можна вважати, що $(\psi_m\circ g_m)(z_0)=0.$ Тоді
з~(\ref{eq2E}) випливає, що
\begin{equation}\label{eq2F}
f_m=\varphi_m\circ g_m= \varphi_m\circ\psi^{\,-1}_m\circ\psi_m\circ
g_m=F_m\circ G_m\,,\quad m=1,2,\ldots\,,
\end{equation}
де $F_m:=\varphi_m\circ\psi^{\,-1}_m,$ $F_m:{\Bbb D}\rightarrow
{\Bbb C},$ і $G_m=\psi_m\circ g_m.$
Очевидно, функція $F_m$ є аналітичною, а $G_m$ -- регулярний
гомеоморфізм класу Соболєва в області $D.$ Зокрема, ${\rm
Im}\,F_m(0)=0$ для всіх $m\in {\Bbb N}.$

\medskip
\textbf{III.} Зауважимо що
\begin{equation}\label{eq7A}
\int\limits_D K_{\mu_{G_m}}(z)\,dm(z) \leqslant \int\limits_D
Q_M(z)\,dm(z) <\infty\,,
\end{equation}
оскільки за умовою $\mu(z)\in M(z)$ при майже всіх $z\in D.$ Отже,
при майже всіх $z\in D$ виконано нерівність
$K_{\mu_{G_m}}(z)\leqslant Q_M(z),$ причому функція $Q_M(z)$ не
залежить від індексу $m=1,2,\ldots $ і інтегровна в $D$ за умовою.

\medskip
\textbf{IV.} Доведемо, що кожне відображення $G_m,$ $m=1,2,\ldots ,$
має неперервне продовження на $\partial D,$ крім того, сім'я
продовжених відображень $\overline{G}_m,$ $m=1,2,\ldots ,$ є
одностайно неперервною в $\overline{D}.$ Дійсно, оскільки
$K_{\mu_{G_m}}(z)\leqslant Q_M(z)$ при майже всіх $z\in D,$
$K_{\mu_{G_m}}\in L^1(D).$ З огляду на це, за~\cite[теорема~3]{KPRS}
(див. також~\cite[теорема~3.1]{LSS}) кожне відображення $G_m,$
$m=1,2,\ldots ,$ є так званим кільцевим $Q_M$-гомеоморфізмом в
$\overline{D}.$ Тоді бажаний висновок є твердженням леми~\ref{lem2}.

\medskip
\textbf{V.} Доведемо також, що обернені гомеоморфізми $G^{\,-1}_m,$
$m=1,2,\ldots ,$ продовжується по неперервності на $\partial {\Bbb
D}$ і сім'я відображень $\{\overline{G}^{\,-1}_m\}_{m=1}^{\infty}$ є
одностайно неперервною в $\overline{\Bbb D}.$  Оскільки за доведеним
у пункті~\textbf{IV} відображення $G_m,$ $m=1,2,\ldots, $ є
кільцевими $K_{\mu_{G_m}}(z)$-гомеоморфізмами в $D,$ обернені до них
відображення $G^{\,-1}_m$ задовольняють співвідношення~(\ref{eq2*A})
(в цьому випадку, $D$ у (\ref{eqA2}) відповідає одиничному кругу
${\Bbb D},$ $f\mapsto G_m,$ $Q\mapsto K_{\mu_{G_m}}(z),$ відповідно,
області $f(D)$ у~(\ref{eq2*A}) відповідає $D$). Оскільки
$G^{\,-1}_m(0)=z_0$ для всіх $m=1,2,\ldots ,$ неперервне продовження
кожного відображення $G^{\,-1}_m$ на $\partial {\Bbb D},$ а також
одностайна неперервність сім'ї відображень
$\{\overline{G}_m^{\,-1}\}_{m=1}^{\infty}$ на $\overline{\Bbb D}$ є
результатом леми~\ref{th4}.

\medskip
\textbf{VI.} Оскільки за доведеним в пункті~\textbf{IV} сім'я
відображень $\{\overline{G}_m\}_{m=1}^{\infty}$ є одностайно
неперервною в~$\overline{D},$ за критерієм Арцела-Асколі існує
зростаюча підпослідовність номерів $m_k,$ $k=1,2,\ldots ,$ така що
послідовність $\overline{G}_{m_k}$ збігається рівномірно в
$\overline{D}$ при $k\rightarrow\infty$ до деякого неперервного
відображення $\overline{G}:\overline{D}\rightarrow \overline{{\Bbb
C}}$ (див., напр., \cite[теорема~20.4]{Va}). Нехай $G\colon
=\overline{G}|_{D}.$ За лемою~\ref{lem1} має місце альтернатива: або
$G$ -- гомеоморфізм області $D$ у ${\Bbb C},$ або $G$ -- стала в
$\overline{\Bbb C}.$ Доведемо, що другий випадок неможливий.
Скористаємося підходом, застосованим при доведенні другої частини
теореми~21.9 в~\cite{Va}. Припустимо супротивне: нехай
$G_{m_k}(x)\rightarrow c=const$ при $k\rightarrow\infty.$ Оскільки
$G_{m_k}(z_0)=0$ при всіх $k=1,2,\ldots ,$ маємо: $c=0.$ З огляду на
пункт~\textbf{V} сім'я відображень $G^{\,-1}_m,$ $m=1,2,\ldots ,$ є
одностайно неперервною в ${\Bbb D}.$ Тоді для довільної точки $z\in
D$
$$h(z, G^{\,-1}_{m_k}(0))=h(G^{\,-1}_{m_k}(G_{m_k}(z)), G^{\,-1}_{m_k}(0))\rightarrow 0$$
при $k\rightarrow\infty,$ що неможливо, бо $z$ -- довільна точка
області~$D.$ Отримана суперечність вказує на те, що $G:D\rightarrow
{\Bbb C}$ -- гомеоморфізм.

\medskip
\textbf{VII.} За доведеним у пункті~\textbf{V} сім'я відображень
$\{\overline{G}_m^{\,-1}\}_{m=1}^{\infty}$ є одностайно неперервною
в $\overline{\Bbb D}.$ Отже, за критерієм Арцела-Асколі (див.,
напр., \cite[теорема~20.4]{Va}) ми також можемо вважати, що
послідовність $\overline{G}^{\,-1}_{m_k}(y),$ $k=1,2,\ldots, $
збігається рівномірно в $\overline{\Bbb D}$ до декого неперервного
відображення $\widetilde{F}:\overline{{\Bbb D}}\rightarrow
\overline{D}$ при $k\rightarrow\infty.$ Встановимо, що
$\widetilde{F}=\overline{G}^{\,-1}.$ Для цього покажемо, що
$G(D)={\Bbb D}.$ Зафіксуємо $y\in {\Bbb D}.$ Оскільки
$G_{m_k}(D)={\Bbb D}$ при всіх $k=1,2,\ldots, $ ми маємо
$G_{m_k}(x_k)=y$ при деякому $x_k\in D.$ Оскільки область $D$
обмежена, можна вважати, що $x_k\rightarrow x_0\in \overline{D}$ при
$k\rightarrow\infty.$ Далі, використовуючи нерівність трикутника і з
огляду на одностайну неперервність
$\{\overline{G}_m\}_{m=1}^{\infty}$ в $\overline{D}$
(пункт~\textbf{IV}), будемо мати:
$$|\overline{G}
(x_0)-y|=|\overline{G}(x_0)-\overline{G}_{m_k}(x_k)|\leqslant
|\overline{G}(x_0)-\overline{G}_{m_k}(x_0)|+|\overline{G}_{m_k}(x_0)
-\overline{G}_{m_k}(x_k)|\rightarrow 0$$
при $k\rightarrow\infty.$ Звідси $\overline{G}(x_0)=y.$ Зауважимо,
що $x_0\in D,$ оскільки $G$ -- гомеоморфізм. В силу довільності
точки $y\in {\Bbb D}$ рівність $G(D)={\Bbb D}$ доведено. В такому
випадку, $G^{\,-1}_{m_k}\rightarrow G^{\,-1}$ локально рівномірно в
${\Bbb D}$ при $k\rightarrow\infty$ (див., напр.,
\cite[лема~3.1]{RSS}). Таким чином, $\widetilde{F}(y)=G^{\,-1}(y)$
при всіх $y\in {\Bbb D}.$ Нарешті, оскільки відображення
$\widetilde{F}$ має неперервне продовження на межу області ${\Bbb
D},$ то в силу єдиності границі в межових точках маємо також
$\widetilde{F}(y)=\overline{G}^{\,-1}(y)$ при всіх $y\in
\overline{{\Bbb D}}.$ Отже, ми довели, що
$\overline{G}^{\,-1}_{m_k}\rightarrow \overline{G}^{\,-1}$
рівномірно в~$\overline{\Bbb D}$ при $k\rightarrow\infty.$

\medskip
\textbf{VIII.} За пунктом~\textbf{VII,} для $y=e^{i\theta}\in
\partial {\Bbb D}$ при $k\rightarrow\infty$ будемо мати:
\begin{equation}\label{eq4E}
{\rm
Re\,}F_{m_k}(e^{i\theta})=\varphi(\overline{G}^{\,-1}_{m_k}(e^{i\theta}))\rightarrow
\varphi(\overline{G}^{\,-1}(e^{i\theta}))
\end{equation}
рівномірно по $\theta\in [0, 2\pi).$ Оскільки за побудовою ${\rm
Im\,}F_{m_k}(0)=0$ при всіх $k=1,2,\ldots,$ за формулою Шварца
(див., напр., \cite[$\S\,$8, гл.~III, частина~3]{GK}) аналітична
функція $F_{m_k}$ однозначно відновлюється по своїй дійсній частині,
а саме,
\begin{equation}\label{eq4A}
F_{m_k}(y)=\frac{1}{2\pi i}\int\limits_{S(0,
1)}\varphi(\overline{G}^{\,-1}_{m_k}(t))\frac{t+y}{t-y}\cdot\frac{dt}{t}\,.
\end{equation}
Покладемо
\begin{equation}\label{eq5A}
F(y):=\frac{1}{2\pi i}\int\limits_{S(0,
1)}\varphi(\overline{G}^{\,-1}(t))\frac{t+y}{t-y}\cdot\frac{dt}{t}\,.
\end{equation}
Нехай $K\subset {\Bbb D}$ -- довільний компакт. З огляду на
співвідношення~(\ref{eq4A}) і~(\ref{eq5A}) ми отримаємо, що для
$z\in K$
\begin{equation}\label{eq11A}
|F_{m_k}(y)-F(y)|\leqslant \frac{1}{2\pi}\int\limits_{S(0,
1)}|\varphi(\overline{G}^{\,-1}_{m_k}(t))-\varphi(
\overline{G}^{\,-1}(t))|\left|\frac{t+y}{t-y}\right|\,|dt|\,.
\end{equation}
Оскільки $K$ -- компакт, знайдеться $0<R_0=R_0(K)<0$ таке, що
$K\subset B(0, R_0).$ Тоді за нерівністю трикутника $|t+y|\leqslant
1+R_0$ і $|t-y|\geqslant |t|-|y|\geqslant 1-R_0$ для всіх $y\in K$ і
всіх $t\in {\Bbb S}^1.$ Тоді
\begin{equation}\label{eq12D}
\left|\frac{t+y}{t-y}\right|\leqslant \frac{1+R_0}{1-R_0}:=M=M(K)\,.
\end{equation}
Зафіксуємо довільне $\varepsilon>0.$ З огляду на умову~(\ref{eq4E})
для числа $\varepsilon^{\,\prime}:=\frac{\varepsilon}{M}$ знайдеться
номер $N=N(\varepsilon, K)\in {\Bbb N}$ такий, що
$|\varphi(\overline{G}^{\,-1}_{m_k}(t))-\varphi(\overline{G}^{\,-1}(t))|<\varepsilon^{\,\prime}$
для всіх $k\geqslant N(\varepsilon).$ Тоді з~(\ref{eq11A})
і~(\ref{eq12D}) випливає, що
\begin{equation}\label{eq13A}
|F_{m_k}(y)-F(y)|<\varepsilon \quad \forall\,\,k\geqslant N\,.
\end{equation}
З нерівності~(\ref{eq13A}) випливає, що послідовність $F_{m_k}$
збігається до функції $F$ локально рівномірно в одиничному крузі.
Зокрема, маємо: ${\rm Im\,}F(0)=0,$ тоді також ${\rm Im\,}f(z_0)=0.$
Зауважимо, що $F$ є аналітичною функцією в ${\Bbb D}$ (див.
зауваження, зроблені в кінці параграфу~8 частини~3 у~\cite{GK}),
причому для $z=re^{i\psi}$
$${\rm Re}\,F(re^{i\psi})=\frac{1}{2\pi}\int\limits_0^{2\pi}
\varphi(\overline{G}^{\,-1}(e^{i\theta}))\frac{1-r^2}{1-2r\cos(\theta-\psi)+r^2}\,d\theta\,.$$
За~\cite[теорема~2, $\S\,$10, гл.~III, частина~3]{GK}
\begin{equation}\label{eq15A}
\lim\limits_{\zeta\rightarrow z}{\rm
Re}\,F(\zeta)=\varphi(\overline{G}^{\,-1}(z))\quad\forall\,\,z\in\partial
{\Bbb D}\,.
\end{equation}
Зауважимо, що функція $F$ є або сталою, або відкрита і дискретна
(див., напр., \cite[гл.~V, розд.~I, пункт~6 і розд.~II,
пункт~5]{St}). Отже, послідовність $f_{m_k}=F_{m_k}\circ G_{m_k}$
збігається локально рівномірно до функції $f=F\circ G,$ яка є
відкритою і дискретною, або сталою функцією, причому, з огляду
на~(\ref{eq15A})
$${\rm Re\,}f(z)={\rm Re\,}F(\overline{G}(z))=\varphi(\overline{G}^{\,-1}(\overline{G}(z)))=\varphi(z)\,.$$
\textbf{IX.} Оскільки за доведеним у пункті~\textbf{VI} відображення
$G$ є гомеоморфізмом, з огляду на \cite[теорема~1]{L$_2$} $G$ є
регулярним розв'язком рівняння~(\ref{eq2C}) з деякою
функцією~$\mu:D\rightarrow {\Bbb D}.$ Оскільки множина точок функції
$F,$ де її якобіан дорівнює нулю, може складатися тільки з
ізольованих точок (див.~\cite[пункти 5 и 6 (II), гл.~V]{St}), у
випадку $F\not\equiv const$ відображення $f$ є регулярним.
Залишилось довести, що $\mu\in \frak{M}_M.$ Якщо $f(z)=c=const$ в
області $D,$ то завдяки умові~(\ref{eq1A}) це є можливим лише у
випадку, коли $f_n(z)=c$ у $D,$ а $\mu_n(z)=0\in M(z)$ при майже
всіх $z\in D.$ В цій ситуації також $\mu(z)=0$ при майже всіх $z\in
D,$ зокрема, $\mu\in \frak{M}_M.$

\medskip
Нехай тепер $f(z)\ne const.$ За доведеним вище відображення $f$ є
регулярним. Оскільки $f_n(z)$ збігається до $f(z)$ локально
рівномірно в області $D$ і, крім того, $f$ має майже всюди відмінний
від нуля якобіан, то за~\cite[лема~1]{L$_2$} $\mu(z)\in {\rm
inv\,co} M_0(z)$ для майже всіх $z\in D,$ де ${\rm inv\,co}\,A$
позначає інваріантно опуклу оболонку множини $A\subset {\Bbb C}$
(див., напр., \cite{Ryaz}), а $M_0(z)$ позначає множину точок
скупчення послідовності $\mu_n(z),$ $n=1,2,\ldots .$ Очевидно, існує
множина $D_0\subset D$ така, що $\mu_n(z)\in M(z)$ і $\mu(z)\in {\rm
inv\,co}\, M_0(z)$ при всіх $z\in D_0$ і всіх $n\in {\Bbb N},$ де
$m(D\setminus D_0)=0.$ Зафіксуємо $z_0\in D_0.$ Нехай $w_0\in
M_0(z_0).$ Тоді існує підпослідовність номерів $n_{k},$
$k=1,2,\ldots ,$ для якої $\mu_{n_k}(z_0)$ є збіжною при
$k\rightarrow\infty$ і $\lim\limits_{k\rightarrow
\infty}\mu_{n_k}(z_0)=w_0.$ Оскільки за припущенням
$\mu_{n_k}(z_0)\in M(z_0)$ при всіх $k=1,2,\ldots ,$ крім того,
множина $M(z_0)$ є замкненою, то $w_0\in M(z_0).$ Отже,
\begin{equation}\label{eq1D}
M_0(z_0)\subset M(z_0)\,.
\end{equation}
Зі співвідношення~(\ref{eq1D}) випливає, що
\begin{equation}\label{eq2J}
{\rm inv\,co}\,M_0(z_0)\subset M(z_0)\,,
\end{equation}
оскільки множинна $M(z_0)$ припускалася інваріантно опуклою. Отже,
$$\mu(z_0)\in {\rm inv\,co}\,M_0(z_0)\subset M(z_0)$$ при майже всіх
$z_0\in D,$ що і потрібно було довести.~$\Box$

\medskip
Для фіксованої функції $Q:{\Bbb C}\rightarrow [0, \infty],$
$Q(z)\equiv 0$ при $z\in {\Bbb C}\setminus D,$ точки $z_0\in D$ і
неперервної функції $\varphi:\partial D\rightarrow {\Bbb R}$
позначимо через $\frak{F}_{\varphi, Q, z_0}(D)$ клас усіх регулярних
розв'язків $f:D\rightarrow{\Bbb C}$ задачі
Діріхле~(\ref{eq2C})--(\ref{eq1A}), які задовольняють умову ${\rm
Im}\,f(z_0)=0$ таких, що $K_{\mu_f}(z)\leqslant Q(z)$ для майже всіх
$z\in D.$ Наступне твердження також можна розглядати як
узагальнення~\cite[теорема~2]{Dyb} на випадок однозв'язних
жорданових областей.

\medskip
\begin{corollary}\label{cor1}
{\sl Нехай $D$ -- деяка однозв'язна жорданова область у ${\Bbb C},$
і нехай функція $\varphi$ у~(\ref{eq1A}) неперервна. Припустимо, що
функція $Q$ є інтегровною в $\overline{D}$ і задовольняє принаймні
одну з умов: або $Q\in FMO(\overline{D}),$ або для кожного $z_0\in
\overline{D}$ існує $\delta_0=\delta(z_0)>0$ таке, що виконано
умову~(\ref{eq2G}), де
$q_{z_0}(t)=\frac{1}{2\pi}\int\limits_0^{2\pi}Q(z_0+e^{i\theta})\,d\theta.$
Тоді сім'я відображень $\frak{F}_{\varphi, Q, z_0}(D)$ є компактною
в $D.$}
\end{corollary}

\medskip
\begin{proof}
Дійсно, в силу рівності
$K_{\mu_f}(z)=\frac{1+|\mu_f(z)|}{1-|\mu_f(z)|}$ умова
$K_{\mu_f}(z)\leqslant Q(z)$ еквівалентна умові $|\mu_f(z)|\leqslant
\frac{Q(z)-1}{Q(z)+1}.$ Тоді $\mu_f(z)\in \overline{B\left(0,
\frac{Q(z)-1}{Q(z)+1}\right)},$ причому множини
$M(z):=\overline{B\left(0, \frac{Q(z)-1}{Q(z)+1}\right)}$ є
замкненими і інваріантно опуклими. В такому випадку, бажане
твердження випливає з теореми~\ref{th2A}.~$\Box$
\end{proof}


КОНТАКТНА ІНФОРМАЦІЯ

\medskip
\noindent{{\bf Олександр Петрович Довгопятий} \\
Житомирський державний університет ім.\ І.~Франко\\
кафедра математичного аналізу, вул. Велика Бердичівська, 40 \\
м.~Житомир, Україна, 10 008 \\
e-mail: alexdov1111111@gmail.com}

\medskip
\noindent{{\bf Євген Олександрович Севостьянов} \\
{\bf 1.} Житомирський державний університет ім.\ І.~Франко\\
кафедра математичного аналізу, вул. Велика Бердичівська, 40 \\
м.~Житомир, Україна, 10 008 \\
{\bf 2.} Інститут прикладної математики і механіки
НАН України, \\
вул.~Добровольського, 1 \\
м.~Слов'янськ, Україна, 84 100\\
e-mail: esevostyanov2009@gmail.com}

\end{document}